\documentclass[11pt,a4paper]{article}
\usepackage{amsmath,amssymb,amsthm,url}
\usepackage{amsfonts,enumerate}
\usepackage[usenames,dvipsnames]{color}
\newcommand\red[1]{#1}
\input{xy}\xyoption{all}
\usepackage[a4paper,left=3cm,right=2cm,top=2cm,bottom=2.8cm]{geometry}

 \def\R{{\mathbb R}}  
\long\def\comment#1\endcomment{}

\def\st{\mathop{\fam0 st}}

\theoremstyle{plain}
\newtheorem{Theorem}{Theorem}
\newtheorem{Conjecture}[Theorem]{Conjecture}

\theoremstyle{definition}
\newtheorem{Remark}[Theorem]{Remark}

\begin{document}
\title{On the metastable Mabillard-Wagner conjecture
\thanks{Research supported by the Russian Foundation for Basic Research Grant No. 15-01-06302, by Simons-IUM Fellowship and by the D. Zimin's Dynasty Foundation Grant.
\newline
I would like to thank I. Mabillard U. and Wagner for helpful discussions.}
}

\def\istsymb{\textrm{a}}
\def\moscowsymb{\textrm{b}}

\author{Arkadiy Skopenkov$^{\moscowsymb}$}

\date{}

\maketitle

{\renewcommand\thefootnote{\moscowsymb}
\footnotetext{Moscow Institute of Physics and Technology,
 and Independent University of Moscow.
Email: \texttt{skopenko@mccme.ru}}
}

\begin{abstract}
The purpose of this note is to attract attention to the following conjecture (metastable $r$-fold Whitney trick) 
by clarifying its status as not having a complete proof, in the sense described in the paper.
{\it Assume that $D=D_1\sqcup\ldots\sqcup D_r$ is disjoint union of $r$ disks of dimension $s$, 
\ $f:D\to B^d$ a proper PL map such that $f\partial D_1\cap\ldots\cap f\partial D_r=\emptyset$,  
\ $rd\ge (r+1)s+3$ and $d\ge s+3$.
If the map
$$f^r:\partial(D_1\times\ldots\times D_r)\to (B^d)^r-\{(x,x,\ldots,x)\in(B^d)^r\ |\ x\in B^d\}$$
extends to $D_1\times\ldots\times D_r$, then there is a PL map $\overline f:D\to B^d$ such that
$$\overline f=f \quad\text{on}\quad D_r\cup\partial D\quad\text{and}\quad
\overline fD_1\cap\ldots\cap \overline fD_r=\emptyset.$$}
\end{abstract}

The purpose of this note is to attract attention to the following conjectures by clarifying their status as 
not having {\it complete} proofs, in the sense described below.
 
Let $B^{d}:=[0,1]^d$ denote the standard PL (piecewise linear) ball and $S^{d-1}=\partial B^d$ the standard PL sphere.
We need to speak about PL balls of different dimensions and we will use the word `disk' for lower-dimensional
objects and `ball' for higher-dimensional ones in order to clarify the distinction (even though, formally, the disk $D^d$ is the same as the ball $B^d$).
%We denote by $\partial M$, respectively $\Int M$, the boundary, respectively the interior, of a manifold $M$.
A map $f\colon M \to B^d$ from a manifold with boundary to a ball is called {\it proper}, 
if $f^{-1}S^{d-1}=\partial M$.
%In this paper we work in the PL category \cite{RS72}, in particular, all disks and balls are PL.
 
%Denote
%$$\diag\phantom{}_r:=\{(x,x,\ldots,x)\in(\R^d)^r\ |\ x\in\R^d\}.$$
%Recall that this map is a $\Sigma_r$-equivariant homotopy equivalence.
%$(\R^d)^r-\diag\sim_{\Sigma_r} S^{d(r-1)-1}_{\Sigma_r}$.

\begin{Conjecture}[Metastable Local Disjunction]\label{l:ldm}
Assume that

$\bullet$ $D=D^{s_1}\sqcup\ldots\sqcup D^{s_r}$ is disjoint union of $r$ disks; 

$\bullet$ $f:D\to B^d$ a proper PL map such that $f\partial D^{s_1}\cap\ldots\cap f\partial D^{s_r}=\emptyset$;

$\bullet$ $rd\ge (s_1+s_2+\ldots+s_r)+s_i+3$ and $d\ge s_i+3$ for each $i$.

If the map
$$f^r:\partial(D^{s_1}\times\ldots\times D^{s_r})\to (B^d)^r-\{(x,x,\ldots,x)\in(B^d)^r\ |\ x\in B^d\}$$
extends to a continuous map of $D^{s_1}\times\ldots\times D^{s_r}$, 
then there is a PL map $\overline f:D\to B^d$ such that
$$\overline f=f \quad\text{on}\quad D^{s_r}\cup\partial D\quad\text{and}\quad
\overline fD^{s_1}\cap\ldots\cap \overline fD^{s_r}=\emptyset.$$
\end{Conjecture}

%Throughout this paper, let $K$ be a finite simplicial complex.
A continuous (or PL) map  $f\colon K\to \R^d$ of a finite simplicial complex is an {\it almost $r$-embedding} if
$f(\sigma_1)\cap \ldots \cap f(\sigma_r)=\emptyset$ whenever $\sigma_1,\ldots,\sigma_r$ are
pairwise disjoint simplices of $K$.

Denote by $\Sigma_r$ the permutation group of $r$ elements.
The group $\Sigma_r$ acts on the set of real $d\times r$-matrices by permuting the columns.
Denote by $S^{d(r-1)-1}_{\Sigma_r}$ the set of the set of real $d\times r$-matrices such that the sum in each row
is zero, and the sum of squares of the matrix elements is 1.
This set is homeomorphic to the sphere of dimension $d(r-1)-1$.

\begin{Conjecture}[Metastable Mabillard-Wagner Conjecture]\label{t:mmw}
Let $s,d,r\ge2$ be integers satisfying $rd\ge(r+1)s +3$ and $K$ a finite $s$-dimensional simplicial complex.
There exists an almost r-embedding $f : K\to\R^d$ if and only if there exists a continuous $\Sigma_r$-equivariant map
%$$K^{\times r}_{\Delta}:=
$$\bigcup \{ \sigma_1 \times \cdots \times \sigma_r
\ : \sigma_i \textrm{ a simplex of }K,\ \sigma_i \cap \sigma_j = \emptyset \mbox{ for every $i \neq j$} \}
\to S^{d(r-1)-1}_{\Sigma_r}.$$
\end{Conjecture}

These conjectures have some interesting corollaries \cite{MW16, MW'}.

These conjectures were claimed as theorems in \cite{MW16, MW'}.
However, in view Remark \ref{r:sp} below I find the proof of the main results in \cite{MW'} {\it incomplete} in the following sense. 
 (Version 1 of \cite{MW'} and the published paper \cite{MW16} contains even less details.) 
In this note I call a proof {\it incomplete} if one mathematician should be able to expect from another

(1) to wait for another (`complete') proof {\it before} using results having such a proof;
%(i.e. whose proof in delicate places contains formulations at the
%level of accuracy described in the bullet points of (c));

(2) to recommend, as a referee, a revision (based on specific comments) {\it before} recommending publication of  results having such a proof;

(3) to work more on such a proof (in particular, send the text privately to a small number of mathematicians working on related problems), {\it before} submitting the text to a refereed journal or to arxiv.% 
\footnote{Unfortunately, shorter formulations of the conclusion of Remark \ref{r:sp} were found to be potentially misleading.}

%We??? find the formulation of this practical statement, with which we explicitly disagree???,
%an important achievement comparable to forming a unanimous opinion.
%on the subject???.
%we could not formulate such statements in a shorter way.
%(because `a gap' or `not up to the standards of refereed journals' were found vague and impolite).

No other meaning of `{\it incomplete}' is meant here.
%In spite of the above criticism,
In particular, I have nothing against publication of {\it conjectures} with {\it incomplete} proofs. 
And I do not mean that Conjectures \ref{l:ldm} and \ref{t:mmw} are wrong, or that the proof of \cite{MW'} cannot be recovered.
%By a {\it proof} I mean proof of a result stated as a theorem, not that it is improper to publish results having
In \cite{Sk17} I give a different short proof (based on working in the smooth category),  
%instead of correcting drawbacks of Remark \ref{r:sp}.b), 
and call Conjecture \ref{t:mmw} Metastable Mabillard-Wagner Theorem.
This would hopefully allow to concentrate not on priority but on mathematics and on the following question: 

{\it which level of accuracy is required to recognize a mathematical proof as complete, in the above sense?}

It is important to answer this question in specific situations, 
so that our journal or arxiv publications will not be treated with suspicion.
 
U. Wagner and I find open publication of this criticism of \cite{MW16, MW'} important to stimulate appearance of a complete proof.
I also find this important to stimulate private discussions (among mathematicians working in the same area) before making claims in arxiv or journal submissions.
This note might also be interesting as an example of an open discussion of a potentially controversial question, 
carried in full mutual respect of participants of the discussion. 

%\red{In case of disagreement this sentence could be moved to (e)
%should be given in a longer form:
%A. Skopenkov finds open publication of the above criticism important to stimulate appearance of a complete proof.
%Indeed, the problem of accurate use of block bundles theory was known to I. Mabillard and U Wagner before
%the submission of version 1 of \cite{MW'} on 5 Jan 2016 but was neither recognized nor completely resolved in
%this or subsequent publications [more history].}
%, if required].}

We discussed the criticism with I. Mabillard and U. Wagner.
I asked them if they agree that {\it the proof of the main results in \cite{MW'} is} incomplete {\it in the above sense.}
I received no answer to this question.
However, I received Remark \ref{r:uw} below. 
\red{It first very politely agrees with all Remark \ref{r:sp}.b but the last two bullet points. 
%(of which the last point is hopefully not significant, and was hidden in \cite{Sk17}, 
%which obstructed its discussion).
I am afraid the rest of Remark \ref{r:uw} gives no information relevant to discussing {\it incompleteness} of the proof. 
However, since Uli Wagner looked at most of my comments in footnotes and suggested to publish Remark \ref{r:uw},
I am glad to include it here.}
 
\begin{Remark}\label{r:sp}
(a) Besides the main result of \cite{MW'}, I do not know any result whose formulation does not not use block bundles theory \cite{RS} but which could be proved using this theory, and whose proof not using this theory is unknown or hard.
(Although Rourke and Sanderson might have been aware of such results while they wrote \cite{RS}.)%  
\footnote{I also have to study a reference to such a result sent to me by I. Mabillard.}
For me it was always hard to apply block bundles theory for such a result, and it was easier to use other means,
even if the result was invented by guessing that in a given situation PL manifolds behave analogously to smooth
ones, cf. \cite[\S4]{Sk02}, \cite[Remark 21]{MW'}, \cite{Sk17}, \cite{CS}.
So application of block bundles theory is a delicate part of the proof.
The corresponding Lemmas 16, 19 and 24 \cite{MW'} are not rigorously stated and proved, see specific remarks 
in (b) below.
These lemmas are important for  the proof.
Accurate statements of the lemmas would be more technical.

Checking proofs and use of the lemmas in their accurate statements would become an important task
for the authors.
As always, in performing this task new problems might or might not be discovered 
\red{(an example of such a discovery is the last but one bullet point of (b) concerning Observation 18).} 
There is no way to see how important or negligible these problems are, except authors doing this work and looking critically on the new text.

(b) $\bullet$ P. 25, proof of Lemma 10, First Part.  Theorem 7 is not applicable because Proposition 13 does not
assert that $\sigma_i\cap\sigma_r$ is a PL manifold.
This assumption is tacitly used in the proof of Proposition 13 (i.e. of Lemma 14, \S4.1) but never checked.
This assumption is non-trivial to check, cf. footnote 11 in version 2 of \cite{AMSW}.

$\bullet$ In Lemma 16 {\it unknottedness}, i.e. the existence of a homeomorphism (or an isotopy), is tacitly replaced by {\it equality}; we know it is dangerous to identify isomorphic objects; cf. use of `$\cong$' not `$=$' in the analogous situations in Lemma 19(2) and in Lemma 20;

$\bullet$ In Lemma 16 `we can assume' is way too informal for this delicate part of the proof;

$\bullet$ In Lemma 16 it is not clear where the restatement of the transversality (the second `i.e.') ends: at `in $B^d$' or at the end of the display formula;

$\bullet$ In the proof of Lemma 16 `follows by Theorem 66' is not clear because no block bundle (required for application of Theorem 66) is given in the statement of Lemma 16;

$\bullet$ In Lemmas 19, 20 and 24 `In the situation given by Lemma NN' is unclear because it is not indicated whether the assumption or the conclusion of Lemma NN is meant;

$\bullet$ In the proof of Lemma 19 `the first property follows from Theorem 70 $<$presumably meaning Proposition 70$>$ from Appendix A' is not clear because Proposition 70 does not assert any triviality (or the existence of any homeomorphism, which the triviality means).

$\bullet$ Observation 18 should be turned into a formal statement having formal proof. 
Otherwise in Observation 18 `Since $k\le s-(r-1)d$' is not clear because
it is not written what are the hypothesis of Observation 18, and in the first line of the proof of Lemma 19  
`follows by ... the above observation' is not clear because it is not written which exactly statement is meant by  
`the above observation'. 
This done, one would see that one needs to prove that `if $2k+1\le s_i+s_r-d$, then the unstated conclusion of Observation 18 holds' (currently this part of the proof is ignored). 
Then one would see that one needs the restriction $2k+2\le s_i+s_r-d$, but  not $2k+1\le s_i+s_r-d$, unless one wants to show that methods of \cite{KM}, not only of the cited paper \cite{Mi}, are applicable.

$\bullet$ The statement 8.3 meaning $B^d\cap fK=B^d\cap f(\sigma_1\sqcup\ldots\sqcup\sigma_r)$ is wrong and
should be replaced by $B^d\cap fK=B^d\cap f\st(\sigma_1\sqcup\ldots\sqcup\sigma_r)$; corresponding changes are required in applications of Lemma 8.
\end{Remark}

%The structure of the proof is at some places not clear.
%!At the end of \S4.2 `proving Lemma 33, this is what the rest of this section is devoted to' is wrong, cf. \S4.5.

\begin{Remark}\label{r:uw} {\it This is U. Wagner's public response to the above criticism, only footnotes are mine.} 

Dear Arkadiy,

%As I told you during our Skype conversation, on January 29, 
I agree with many of the specific criticisms (all but the last two bullet points in Remark 3(b)).  
%both the older ones 
%that you mentioned in our Skype conversation on January 29 and the new ones;  
%that you pointed out since then; 
I agree that these statements and formulations should be made more precise. 
(I am not sure regarding the modified dimension restriction you mention in the last
but one bullet points in Remark 3(b), I have to think about this more.%
\footnote{AS: This corresponds to a previous version of the present note where my remark was stated
without the last line mentioning  references \cite{KM, Mi}.}
)

I think these are helpful and valid criticisms, and we will address them in the next revision. 
However, as I tried to explain by skype, I think these things can be easily fixed (in the sense that 
no new ideas are needed)  to arrive at a version of the proof that is hopefully
complete according to your definition.% 
\footnote{\label{f:unc} AS: This sentence and the next paragraph are misleading because they convey the author's disagreement with the statement `the proof of \cite{MW'} is {\it incomplete} (in the practical sense described above after Conjecture \ref{t:mmw})',  without explicitly stating disagreement  and so taking responsibility for such a statement.
Instead of explicit disagreement with the {\it incomplete}, the author introduces vague notions of `easily fixed in the sense that no new ideas are needed' and `gaps as opposed to issues of improving the presentation'.
\newline
E.g. working in the smooth category as in \cite{Sk17} may be regarded as not a new idea (but rather as
`un-introduction' of the new idea of making an argument, originally introduced in the smooth category, work in the PL category). However, the author would hardly consider \cite{Sk17} as an `easy fix'.
We also know that {\it realization} of an idea may be harder than its {\it introduction}.
We can learn that there are no `gaps' in a new proof only if the proof has `presentation', which is not
{\it incomplete} (in the above sense), cf. the second paragraph of Remark \ref{r:sp}.a.}
A number of your remarks can be addressed by straightforward rewording or changing of punctuation.% 
\footnote{\label{f:exp} AS: I welcome these explanations.
However, since the author neither privately distributed the corresponding update of \cite{MW'}, nor submitted it to arxiv, nor gave a list of {\it all} required changes, nor explicitly disagreed with the {\it incompleteness} of the proof, these explanations are not relevant to discussing the {\it incompleteness}. 
So I suggested to omit this part of discussion until the revised version will be ready. 
\newline
I would also call these explanations `continuing to write a complete proof' rather than `rewording or changing of punctuation'.
Cf. the second paragraph of Remark \ref{r:sp}.a.
}
%Thus I solidarize with the author's earlier public statement
%`{\it it will be most productive to discuss that revised version of the proof when it is ready}'.}
In his email from yesterday (February 13), Isaac also gave quite detailed written explanations concerning your remarks; in the following, I repeat some of these explanations, with Isaac's permission:

\begin{itemize}
\item \emph{in Lemmas 19, 20 and 24 ``In the situation given by Lemma NN'' is unclear because
it is not indicated whether the assumption or the conclusion of Lemma NN is meant}

It means under the same hypotheses and using the same notation as in Lemma 10.

\item \emph{in Lemma 16 it is not clear where the restatement of the transversality (the second `i.e.') ends:
at `in $B^d$' or at the end of the display formula;}

it starts at ``i.e.'' and ends at the end of the sentence.

\item \emph{Theorem 7 is not applicable because Proposition 13 does not assert that $\sigma_i\cap \sigma_j$ is a PL manifold.}

$\sigma_i$ and $\sigma_j$ are PL-balls properly embedded inside of a bigger ball, and by general position it is a PL manifold See, e.g., Theorem 1 in Armstrong \& Zeeman, Transversality for PL Manifolds.
\footnote{AS: In the same-day e-mail I recalled the rest of my remark: `this assumption is non-trivial to check, cf. footnote 11 in version 2 of \cite{AMSW}'.}

\item \emph{in Lemma 16 `we can assume' is way too informal for this delicate part of the proof;}

It means ``after an epsilon-perturbation''

\item \emph{in the proof of Lemma 16 `follows by Theorem 66' is not clear because no block bundle
(required for application of Theorem 66) is given in the statement of Lemma 16.}

The block bundle in question is $\sigma^r \times \varepsilon [-1,1]^{d-s_r}$,
which exists by unknottedness in codimension 3.

\item \emph{in the proof of Lemma 19 `the first property follows from Theorem 70 (presumably meaning Proposition 70) from Appendix A' is not clear because Proposition 70 does not assert any triviality (or the existence of any homeomorphism, which the triviality means).}

The last sentence of Proposition 70 is ``Then any normal bundle over $S^k$ in $M$ is trivial''. This means that any normal (block) bundle over $S^k$ in $M$ (in the sense of Theorem 54) is trivial in the sense of Remark 48.(b) and Def 49. (i.e., the existence of a homeomorphism).
\end{itemize}

%\item \emph{These lemmas are used many times throughout the proof.}
%None of these lemmas are used more than once. They are parts of a bigger proof that was eventually
%split into several pieces. They are assemble sequentially once to build Lemma 14.(a)

For these reasons, I am not convinced that your comments constitute ``gaps'' in our proof
(as opposed to issues of improving the presentation) or that they justify your statement
that our theorems should be considered unproved conjectures. 

I understand that you maybe do not consider these explanations sufficient;%
\footnote{AS: See footnote \ref{f:exp}.} 
or maybe your position is (as I believe you said during our skype meeting) that this amounts to
recovering a complete proof, rather than evidence that your comments do not constitute gaps.%
\footnote{AS: I did not said/wrote this and I do not understand this explanation of my position.  
Sorry for my insufficient knowledge of English.}

From our conversations, I also understand that you consider the words ``easy to fix'', ``gap'',
and ``no new ideas needed'' as vague and not practical. I agree that these it may be hard
to arrive at a general definition of these terms, but I do think that mathematicians often use
these words and reach agreements in specific cases.% 
\footnote{AS: These expressions without their explanations in practical terms are appropriate in a text meant 
to be vague, but are misleading in a seemingly controversial reply to a practical question. 
It is stating a practical question, but not bringing many vague notions into discussion, that helps to reach an agreement (or an explicit disagreement, which is also valuable).}

You suggest a definition of an ``incomplete proof'' in your note.
I am not convinced that your definition is more precise or more practical.%  
\footnote{AS: The above notion of {\it incomplete} is of course not precise. 
However, it is explained in terms of practical decisions which mathematicians have regularly to make.}
It is also based on and uses as defining properties social notions of what one
``mathematician should be able to expect from another'', in particular

\begin{itemize}
\item[(2)] ``to recommend, as a referee (based on specific comments) before
recommending publication of results having such a proof;''
\end{itemize}

I am not sure that this serves as a definition of ``incomplete''.%  
\footnote{AS: It serves by definition.} 
I have certainly both written and received referee reports that recommended revisions without anybody stating 
that this meant the proofs were incomplete or the results unproved.%
\footnote{AS: This is a confusion of a general vague notion of `incomplete' and the practical explanation of 
`incomplete' above. 
Cf. ``No other meaning of `{\it incomplete}' is meant here'' above.
The author can replace the word `incomplete' by `XYZ' in the practical explanation, if it would make it easier to answer the question: {\it do you agree that the proof of the main results in
\cite{MW'} is} XYZ {\it in the above sense?}}

Since we were not able to reach an agreement on this, 
%and you did not want to wait for the revision,% 
%\footnote{AS: It was our joint decision that I should `go ahead and post this note on the arxiv'. 
%In my opinion, I waited too long for a revision having a {\it complete} proof. 
%I cannot give more details because we agreed to omit history from this public discussion.}
I agree that you should make your criticism public, and we work on a revision. 
To me that seems the most productive way forward, and hopefully once the revision is ready we can reach 
an agreement regarding the status of our results and proofs.

Best regards,
Uli
\end{Remark}

\comment

In one point, however, I think I disagree with you:
\footnote{AS: With which exactly statement among written above? }
I see a difference between an arXiv preprint and a final journal version.
Specifically, I think there is value in putting a preprint on the arxiv, even if it is not as polished as one would like, precisely
to get feedback (including from colleagues one does not know who may find the ideas interesting).
\footnote{AS: There is no disagreement because I never stated the opposite.
This remark is misleading because `not as polished as one would like' is vague, see footnote \ref{f:unc},
and because the important practical difference between stating a conjecture and stating a theorem is ignored.
I suggested to replace this paragraph to a statement of whether in view of Remark \ref{r:sp} the author finds
the statement and proof of the main results in \cite{MW'} {\it arxiv-conjectural} (I do).
I call a proof {\it arxiv-conjectural} if one mathematician should be able to expect from another, in a public
arxiv submission, to state a result having such a proof as a conjecture (rather than as a theorem).
No other meaning of `arxiv-conjectural' is meant here.
In particular, I do not mean that at a seminar talk, a result having such a proof should be called a conjecture
(rather than a theorem).}

\endcomment

\end{document}